\documentclass{article}
\usepackage[latin1]{inputenc}
\usepackage{t1enc}
\usepackage{a4}
\usepackage{amsmath}
\usepackage{amssymb}
\usepackage{amsmath,theorem,epsf,psfrag}
\usepackage[pdftex]{graphicx}
\newtheorem{theorem}{Theorem}
\newtheorem{definition}{Definition}

\newtheorem{corollary}{Corollary}

\newtheorem{proposition}{Proposition}

\begin{document}
\begin{center}
{\Large
 Joseph Mecke's last fragmentary manuscripts - a compilation} 
 
 \vspace{1cm} 

{Joseph Mecke} (1938--2014)

\vspace{0.5cm}

{Werner Nagel}\\
{Friedrich-Schiller-Universit\"at Jena,\\
Institut f\"ur Stochastik,\\
Ernst-Abbe-Platz 2,
07743 Jena, Germany.\\
Email: werner.nagel@uni-jena.de}  

\vspace{0.5cm}

{Viola Wei\ss}\\
{Ernst-Abbe-Hochschule Jena},\\
FB Grundlagenwissenschaften.
Carl-Zeiss-Promenade 2,
07745 Jena, Germany.\\
Email: Viola.Weiss@eah-jena.de

\end{center}

\begin{abstract}
Summarizing results from Joseph Mecke's last fragmentary manuscripts, the generating function and the Laplace transform for nonnegative random variables are considered. The concept of thickening of a random variable, as an inverse operation to thinning (which is usually applied to point processes) is introduced, based on generating functions, and a characterization of thickable random variables is given. Further, some new relations between exponential distributions and their interpretation in terms of Poisson point processes are derived with the help of the Laplace transform.
\end{abstract}


\vspace{.5cm}


\section{Introduction}
Joseph Mecke passed away in February 2014, a few days after his 76th birthday. Until his last days, he was dealing with mathematical problems, and he wrote fragments of manuscripts, saved on his computer. His brother, Norbert Mecke, was able to identify the corresponding files; he handed them over to W.N. and V.W., in order to see whether some of the material can be published. The present paper is the result of this compilation.

As emphasized in the introduction of \cite{nw03}, Joseph Mecke preferred to work deep into problems in order to reach a clear insight and a maximum of mathematical elegance. After his last published in a journal \cite{mecke10}, he formulated several new ideas and a wider working agenda. The fragments compiled here date from July 2011 to December 2011 and then from February 2013 to June 2013.

Joseph Mecke made outstanding contributions to the theory of point processes, mainly in the 1960s and early 1970s. Nowadays the Campbell-Mecke formula (Mecke himself referred to it as the 'refined Campbell formula') and the Slivnyak-Mecke formula (see e.g. \cite{sw},  referred to as the Mecke formula in \cite{lastpen}) are cited oftentimes. Since the late 1970-ies, Joseph Mecke worked in stochastic geometry, a field in which he applied the point process theory  strikingly. Thus he contributed to a sound mathematical foundation of this field,  proving rigorously quite a few new results. 

At a first glance, the content of the present paper -- involving  generating functions and the Laplace transform for nonnegative random variables -- seems to be far away from the main subjects of Joseph Mecke's work, described above.
It is not so surprising, however, because 
in his earlier work he applied and appreciated these powerful tools. Although they  appear only occasionally in his published proofs during a long career, this use often gave deeper insight into a problem.

We (W.N. and V.W.) remember a situation in a seminar (in 2006) when we dealt with the length distribution of I-segments in planar STIT tessellations. We had found an expression for the density of this distribution which looked rather strange and we had no clue how to interpret it. Joseph Mecke immediately started his calculation (using the Laplace transform) and soon he revealed this 'mysterious' distribution as a mixture of exponential distributions. Meanwhile, much more is known about STIT tessellations, and there are other methods to prove the mentioned result. But Joseph Mecke opened a door -- as he did it in many other cases.

Probably, the present paper will inspire other mathematicians to study and to generalize some of the problems which Joseph Mecke considered.

\section{Nonnegative integer-valued random variables and (probability) generating functions}

\subsection{Generating function}

We denote ${\mathbb N}_0 =\{0,1,2,\ldots \}$, ${\mathbb N}={\mathbb N}_0\setminus \{0\}$, and ${\bf 1}\{ \cdot \}$ the indicator function which has the value $1$, if the condition in curly brackets is fulfilled and with value $0$ otherwise.

Generating functions are widely  used in mathematics and they play also an important role in probability theory. In this paper they are considered for nonnegative integer-valued  random variables to introduce later the concept of thinning and thickening of a random variable.
Let $\zeta$ be a discrete random variable taking values in ${\mathbb N}_0$ with distribution
\begin{equation}\label{eq:distrzeta}
{\mathfrak L}( \zeta )= \sum_{k=0}^\infty a_k \delta_k ,
\end{equation} 
where $a_k\geq 0$, $\sum_{k=0}^\infty a_k =1$ and $\delta_k$ the Dirac measure assigning mass 1 to $k$.
The corresponding generating function $G :[0,1]\to [0,1]$ is defined by
\begin{equation}\label{eq:nonnegint}
G(x) = \mathbb E (x^{\zeta}) = \sum_{k=0}^\infty a_k \, x^k,\ \quad 0\leq x\leq 1 .
\end{equation}

Note that in the following we will consider the series on the right-hand side also for general $x\in {\mathbb R}$ if it is defined. 

Recall that a function  $G:[0,1]\to [0,1]$ is a probability generating function of a nonnegative integer-valued random variable if and only if $G(1)=1,\ \lim_{x \nearrow 1}G(x) = 1 ,\ G(0) \geq 0$ and all derivatives of $G$ are finite and nonnegative on $[0,1)$ (see e.g. \cite{fristgray}).  Furthermore, the uniqueness theorem conveys that two random variables have the same generating function if and only if they have identical distributions. 

{\bf Examples:}
\begin{enumerate}
\item[(a)]
If $\zeta$ is almost surely (a.s.) constant, $P(\zeta =m)=1$ for some $m\in {\mathbb N}_0$, then $G(x)=x^m $.

\item[(b)] If $\zeta$ has a two-point distribution, ${\mathfrak L}(\zeta)= (1-r)\delta_m + r\delta_n$, $m,n\in {\mathbb N}_0$, $r \in (0,1)$, then the generating function is $G(x)=(1-r)x^m + rx^n$.

\item[(c)] For a random variable $\zeta$ which is Poisson-distributed with parameter $\lambda >0$ we have  $G(x)= {\rm e}^{\lambda (x-1)}$.

\item[(d)] If $\zeta$ has a binomial distribution with parameters $n\in {\mathbb N}$ and $r \in (0,1)$ the generating function is $G(x)=(1-r + rx)^n$, which is the $n$-th power of a generating function of a  Bernoulli random variable with parameter $r$. This follows immediately from the property that the generating function of the sum $\zeta_1 + \zeta_2$ of two independent random variables  is the product of the two generation functions of $\zeta_1$, $\zeta_2$.

\item[(e)] The generating function of a geometric random variable $\zeta$ with parameter $r$ and distribution ${\mathfrak L}(\zeta)=\sum_{k=0}^\infty (1-r)^k r\, \delta_k$ is $G(x)= {\frac{r}{1-x(1-r)}}$.
A negative binomial random variable with distribution 
${\mathfrak L}(\zeta)=\sum_{k=0}^\infty {  -n \choose k}  (r-1)^k r^n \, \delta_k$ (parameters $r \in (0,1)$ and $n \in (0,\infty )$)  has the generating function $ G(x)= \left(\frac{r}{1-x(1-r)} \right) ^n$ which is the $n$-th power of the generating function of a geometric random variable with parameter $r$. 

\end{enumerate}

\hfill $\Box$

\subsection{Thinning and thickening}

Thinning is an operation applied to point processes,  see \cite{dvj} and the references therein. Given a realization, for each single point it is decided (independently of the other points) whether it survives or not. If the survival probability is $p$ for all points, and if $\zeta$ is the (finite) random number of points before thinning, then the distribution of the number of the thinned point process is described in Definition \ref{def:thin}. In this definition, thinning is introduced for arbitrary nonnegative integer-valued random variables. And one can ask whether there is an inverse operation to thinning. So, given a nonnegative integer-valued random variable, can this be the result of thinning of a 'thicker' one, and if so, what is its distribution?  This will be formalized in Definition \ref{def:thickable} and studied in this section.

\begin{definition}\label{def:thin}
Let $\alpha_1, \alpha_2,\ldots $ be independent and identically distributed random variables with the two-point distribution 
${\mathfrak L}(\alpha_k)= (1-p)\delta_0 + p\delta_1$, $k\in {\mathbb N}$.
For a nonnegative integer-valued  random variable $\zeta$ the {\em thinning with parameter $p\in (0,1)$}  is defined as the random variable
\begin{equation}
{\cal D}_p \zeta =\sum_{k=1}^{\zeta} \alpha_k .
\end{equation}
\end{definition}

If $\zeta$ has the distribution given in (\ref{eq:distrzeta}), then the distribution of ${\cal D}_p \zeta$ can be written as
\begin{equation}
{\mathfrak L}({\cal D}_p \zeta ) = \sum_{k=0}^\infty a_k \left( (1-p)\delta_0 + p\delta_1  \right)^{\star k} ,
\end{equation}
where $\star $ denotes the convolution of measures. This means that ${\mathfrak L}({\cal D}_p \zeta )$ is the mixture of binomial distributions with parameters $k$ and $p$, weighted with $a_k$, respectively.

Because $P({\cal D}_p \zeta = m) = \sum_{k=m}^{\infty} a_k { k \choose m }  p^m (1-p)^{k-m}  $, a straightforward calculation yields the generating function $G_p$ of ${\cal D}_p \zeta$,
\begin{equation}\label{eq:genthinning}
G_p(x) =  G(1-p + px).
\end{equation}

{\bf Examples:}
\begin{enumerate}
\item[(a)] For an a.s. constant $\zeta$, $P(\zeta =m)=1$ for some $m\in {\mathbb N}_0$, the generating function of the thinning is $G_p(x)=(1-p+px)^m$.  The uniqueness theorem yields that ${\cal D}_p \zeta$ has a binomial distribution with parameters $m$ and $p$, which is obvious in this case. 
\item[(b)] For $\zeta$ with a two-point distribution, ${\mathfrak L}(\zeta)= (1-r)\delta_m + r\delta_n$,  the thinning ${\cal D}_p \zeta $ is a mixture of two binomial distributions with weights $1-r$ and $r$ and parameters $m$,$n$ respectively, and $p$.

\item[(c)] For a Poisson-distributed $\zeta$ with parameter $\lambda >0$ the generating function of its thinning is
$G_p (x) = {\rm e}^{\lambda p(x-1)}$, which confirms the well-known fact, that ${\cal D}_p \zeta $ has again a Poisson-distribution with parameter 
$p \lambda $.

\item[(d)] If $\zeta$ has a binomial distribution with parameters $n$ and $r$, then the thinning ${\cal D}_p \zeta $ is again binomial distributed but with parameters $n$ and $pr$.

\item[(e)] Also for geometric and negative binomial distributions thinning retaines the type of the distribution. In both cases the parameter $r$ of $\zeta$ changes over to $q=\frac{r}{r +p(1-r)}$ of ${\cal D}_p \zeta $.  
\end{enumerate}

\hfill $\Box$

Now consider $G(x)$, represented by a series as in (\ref{eq:nonnegint}), for arbitrary $x\in {\mathbb R}$ if the value of this series is defined. For $0<p<1$ let us formally modify the function (\ref{eq:genthinning}) to
\begin{equation}\label{eq:genthicking}
G_{\frac{1}{p}} (x) =  G\left( 1-\frac{1}{p} + \frac{1}{p}x\right), \quad \ \ \mbox{ if this is defined for all } 0\leq x \leq 1,
\end{equation}

As a function of $x$, this is not necessarily a generating function of a random variable.

\begin{definition} \label{def:thickable}
Let $\zeta$ be a nonnegative integer-valued random variable with generating function $G$. We say that $\zeta$ is {\em $p$-thickable} for $0<p<1$, if the function $G_{\frac{1}{p}} (x) =  G(1-\frac{1}{p} + \frac{1}{p}x)$ is defined for all $x \in [0,1]$, and if it is the generating function of a nonnegative integer-valued random variable. Such a variable will be denoted by  ${\cal D}_{\frac{1}{p}} \zeta $.
\end{definition} 

A combination of formulas (\ref{eq:genthinning}) and (\ref{eq:genthicking}) yields that

$$
\left(G_{\frac{1}{p}}\right)_p = \left( G_p\right)_{\frac{1}{p}}  =G.
$$
This means that thinning and thickening are somehow mutually inverse operations. But note, that $\left(G_{\frac{1}{p}}\right)_p =G$ is meaningful only for those $p$ for which the distribution is $p$-thickable. The other equation, $ \left( G_p\right)_{\frac{1}{p}}  =G$, holds for all $0<p<1$. This confirms the meaning of thickening as the inverse operation of thinning. 

First we investigate the nonnegative integer-valued  random variables given in the examples above whether they are $p$-thickable or not. 
From the characterization of a generating function it follows that all the (right-hand side) derivatives at $x=0$ are nonnegative. Hence if, for a fixed  $p\in (0,1)$, the function $G_{\frac{1}{p}}$ given in (\ref{eq:genthicking}) is the generating function of a nonnegative integer-valued random variable then
\begin{equation}\label{eq:thickgenerat}
G^{(\ell )}\left( 1-\frac{1}{p}\right) \geq 0 \quad  \mbox{ for all } \ell =0,1,2, \ldots , \
\end{equation}
where $G^{(\ell )}$ denotes the $\ell $-th derivative of $G$. \\[.3cm]

{\bf Examples:}

\begin{enumerate}
\item[(a)] An a.s. constant $\zeta$ with $P(\zeta =m)=1$ is $p$-thickable for all $0<p<1$ if $m=0$, and it is not thickable for any $0<p<1$ if $m$ is a positive integer. For $m=0$ we have $G=G_p=G_{\frac{1}{p}}=1$. In contrast, if $m>0$, then $G(x)=x^m$ and hence $G_{\frac{1}{p}} (x) =  G(1-\frac{1}{p} + \frac{1}{p}x)=(1-\frac{1}{p} + \frac{1}{p}x)^m$. If $m$ is odd, then
$G_{\frac{1}{p}} (x)<0$ for $0\leq x<1-p$. 
And, if $m$ is even, then the first derivative at $x=0$ is negative. This contradicts the necessary condition given in (\ref{eq:thickgenerat}). 

\item[(b)] Analogous considerations show that $\zeta$ with a two-point distribution is not thickable for any $0<p<1$.

\item[(c)] If $\zeta$ is Poisson-distributed with parameter $\lambda >0$ then 
$G_{\frac{1}{p}} (x) = {\rm e}^{\lambda \frac{1}{p}(x-1)}$ which yields that ${\cal D}_{\frac{1}{p}} \zeta $ has a Poisson-distribution with parameter 
$\frac{1}{p} \lambda $. Therefore the Poisson-distributions are $p$-thickable for all $0<p<1$.

\item[(d)] If $\zeta$ has a binomial distribution with parameters $n$ and $r$, i.e. the generating function is $G(x)=(1-r + rx)^n$, then $\zeta$ is $p$-thickable if and only if $r\leq p <1$. For $r>p$ the function $G_{\frac{1}{p}}=(1 - r \frac{1}{p} + r \frac{1}{p}x)^n$ is no longer a generating function. This follows with the same argument for the derivative which was given for $\zeta$ a.s. constant.  For $r\leq p <1$ the $p$-thickening of $\zeta$ is again binomially distributed with parameters $n$ and $\frac{r}{p}$. In particular, the $r$-thickening of $\zeta$ is the constant $n$.

\item[(e)] As in the case of thinning also thickening retains the type of geometric and negative binomial distributions. Thickening is possible for all $p \in (0,1)$ and the new parameter for ${\cal D}_{\frac{1}{p}} \zeta $ is $\frac{r}{r +\frac{1-r}{p}}$.
\end{enumerate}

\subsection{Characterization of unbounded thickability}

In the examples above we have seen that some of the nonnegative integer-valued random variables are $p$-thickable for all $0<p<1$  and others only for some $p$. 

\begin{definition}
A nonnegative integer-valued random variable is called {\em unbounded thickable} if it is $p$-thickable for all $p \in (0,1)$.
\end{definition}

Random variables with a Poisson, a geometrical or a negative binomial distribution  are unbounded thickable. A random variable with a binomial distribution is not unbounded thickable. 
In the following theorem the class of unbounded thickable random variables is described.

\begin{theorem} (Characterization of unbounded thickability)\label{th:unbthick}

A nonnegative integer-valued random variable $\zeta$   with generating function $G$ as in (\ref{eq:nonnegint}) is $p$-thickable {for all} $p\in (0,1)$ if and only if it has a Cox distribution (a mixture of Poisson distributions and the constant 0), i.e. if and only if there exists a probability measure $Q$ on $[0,\infty )$ such that 
\begin{equation}\label{eq:unbthick}
G(x)=\int\limits_{[0,\infty )} {\rm e}^{t(x-1)} Q({\rm d}t), \quad 0\leq x\leq 1.
\end{equation}
\end{theorem}

Proof:

If a function $G$ satisfies (\ref{eq:unbthick}), then obviously $G(1)=1$, $G(0)\geq 0$ and 
$ \lim_{x \nearrow 1}G(x) = 1$, because for all  $0\leq x \leq 1$,  $t\geq 0$, the function ${\rm e}^{t(x-1)}$ is monotone in $x$ and $0< {\rm e}^{t(x-1)} \leq 1$ . Furthermore, because for all $\ell =0,1,\ldots$ the function
$t^\ell {\rm e}^{t(x-1)}$, $t\geq 0$ can be dominated on $(0,\infty )$ by a constant, we obtain that the derivatives $G^{(\ell )}(x) =\int t^\ell {\rm e}^{t(x-1)} Q({\rm d}t)\geq 0$ and they are finite. Hence, $G$ is indeed the generating function of a nonnegative integer-valued random variable. With analogous arguments it can be shown, that also
$G_{\frac{1}{p}} (x) =  G(1-\frac{1}{p} + \frac{1}{p}x)=\int {\rm e}^{\frac{t}{p}(x-1)} Q({\rm d}t)$ is the generating function of a nonnegative integer-valued random variable.

Now we show that (\ref{eq:unbthick}) is necessary for unbounded thickability.
If $\zeta$ is $p$-thickable for all $p\in (0,1)$, then  according to (\ref{eq:thickgenerat}) 
$$
G^{(\ell )}(x) \geq 0 \quad  \mbox{ for all } x <0,\ \ell =0,1,2, \ldots 
$$
For $s\geq 0$ we define
$$
L(s):= G(1-s),
$$
which implies
$$
(-1)^\ell L(s)^{(\ell )}(s)\geq 0 \mbox{ for all } s>0,\ \ell =0,1,2, \ldots
$$
i.e. $L$ is completely monotone on $(0,\infty )$. Furthermore, $L$ is right-continuous at $0$ (because the generating function $G$ is left-continuous at 1) and $L(0)=G(1)=1$. Hence the characterization theorem for Laplace transforms (also referred to as moment generating functions; see e.g. \cite{fristgray}) yields that $L$ is the Laplace transform of a probability measure $Q$ on the half-axis $[0,\infty )$, and hence
$$
G(1-s)=\int {\rm e}^{-ts} Q({\rm d}t), \quad s\geq 0,
$$ 
or, equivalently,
$$
G(x)=\int {\rm e}^{t(x-1)} Q({\rm d}t), \quad x\leq 1.
$$ 
\hfill $\Box$

{\bf Examples:} Referring to the examples above, special Cox distributions are:
\begin{enumerate}
\item[(c)] The Poisson distribution with parameter $\lambda$, and according to (\ref{eq:unbthick}), $Q=\delta_{\lambda}$.

\item[(e)] The negative binomial distribution with parameters $n$ and $r$, where $Q$ is the gamma distribution with parameters $n$ and  $ \frac{r}{1-r} $. In the particular case of a geometric distribution, we have $n=1$ and hence $Q$ is the exponential distribution with parameter $ \frac{r}{1-r} $. 
\end{enumerate}

\subsection{Relations to point processes}

In \cite{amb}, R.V. Ambartzumian  introduced the concept of $1/p$-{\em condensation} of point processes ($p\in (0,1]$) as the inverse operation to thinning. Hence condensation is also related to splitting of point processes. Moreover, he provided a sufficient condition for 2-condensability (which is related to 1/2-thickability considered in the present paper) of point processes in ${\mathbb R}^d$, $d\geq 1$. It remains an open problem to characterize the class of all point processes which are $1/p$-condensable if a value $p\in (0,1)$ is fixed.
Recently,  thinning, splitting and condensation were studied in \cite{NeRaZe, Raf}, and in \cite{Raf} a generalized concept of thinning is introduced, both for nonnegative integer-valued random variables and for point processes. 

In an early paper \cite{mecke68} (where thinning is named 'Ausw{\"u}rfelverfahren'), J. Mecke  already proved  that a point process $\Phi$ on the real axis ${\mathbb R}$ is a Cox process if and only if {\em for any} $p\in (0,1]$ exists a point process $\Phi_p$, such that the $p$-thinning of $\Phi_p$ has the same distribution as $\Phi$ (Satz 4.2 ibidem). This means that a point process on the real axis is unbounded (i.e. for all $p\in (0,1]$) condensable if and only if it is a Cox process.
This result immediately implies  Theorem \ref{th:unbthick} of the present paper. But the proof given here is much shorter and more elegant than that one in \cite{mecke68}. And vice versa, with the help of the generating functional for point processes (see e.g. \cite{meckehabil} or \cite{dvj}), one can easily deduce Satz 4.2 in \cite{mecke68} from Theorem \ref{th:unbthick}.

\subsection{M-transform} \label{sec:M}

The generating function $G$ of a nonnegative integer-valued random variable $\zeta$ in (\ref{eq:nonnegint}) can also be interpreted as the cumulative distribution function (c.d.f.) of a probability measure concentrated on the interval $[0,1]$. Consequently,  in this section we consider the problem how to find for a given $\zeta$ a random variable $\xi$ whose c.d.f. $F_{\xi}$ coincides with $G$ on $[0,1]$. To avoid complications due to $P(\zeta =0)>0$, i.e. $a_0>0$, in this section only positive integer-valued random variables with values in $\mathbb N$ are considered.

\begin{proposition} \label{prop:max}
 Let  $\eta_1,\eta_2,\ldots $ be i.i.d. random variables with uniform distribution on the interval $(0,1)$ and $\zeta$,  independent of this sequence, a positive integer-valued random variable with generating function $G$ as in (\ref{eq:nonnegint}) with $a_0=0$. Then the random variable
$$
\xi :=  \max \{ \eta_1, \ldots ,\eta_\zeta \}
$$
has the c.d.f. $F_\xi$ with $F_\xi (x)=G(x)$ for all $x\in [0,1]$.
\end{proposition}

 Proof: 
 
 Straightforward calculations yield for $0\leq x\leq 1$
$$
\begin{array}{rrl}
F_\xi (x) &=& P(\xi \leq x) \\[2mm]
          &=& \displaystyle{\sum_{k=1}^\infty} P \big( \max \{\eta_1,\ldots , \eta_\zeta \}\leq x | \zeta =k\big) \cdot P(\zeta =k) \\[5mm]
          &=& \displaystyle{\sum_{k=1}^\infty } P\big( \max \{\eta_1,\ldots , \eta_k \} \leq x \big) \cdot P(\zeta =k) \\[5mm]
          &=& \displaystyle{\sum_{k=1}^\infty } x^k \, a_k . 
\end{array}
$$ \hfill $\Box$

Under the assumptions of Proposition (\ref{prop:max}) the following transform of a positive integer-valued random variable can be specified.

\begin{definition}\label{def:Mtrans}
The  random variable $\xi$ on $(0,1)$ given by the transform   $$\xi (\zeta , \eta_1,\eta_2,\ldots )= \max \{ \eta_1, \ldots ,\eta_\zeta \}$$ is called the {\em $M$-transform  of $\zeta$} (maximum-transform), i.e. $\xi = M[\zeta]$.
\end{definition}

{\bf Remark:}
The $M$-transform can also be used in terms of distributions. Let $Q= \sum_{k=1}^\infty a_k \delta_k$, $a_k\geq 0$, $\sum_{k=1}^\infty a_k =1$, be a probability distribution on ${\mathbb N}$ and $U$ the uniform distribution on the interval $(0,1)$. 
Then the $M$-transform of $Q$ is the image measure $M[Q]$ of the product measure $Q\otimes U$ under the transform $(k,t)\mapsto t^{\frac{1}{k}}$, $k\in {\mathbb N}$, $0<t<1$.

Under the assumptions of this definition the transformation theorem for integrals immediately yields
\begin{equation}\label{eq:Mtrans}
\int_{(0,1)} f(s) M[Q]({\rm d}s) =\sum_{k=1}^\infty a_k \int_{(0,1)} f\left( t^{\frac{1}{k}}\right) {\rm d}t \quad\mbox{ for all measurable } f:(0,1)\to [0,1],
\end{equation}
where on the right-hand side ${\rm d}t$ means integration with respect to the Lebesgue measure.

An immediate conclusion from Definition \ref{def:Mtrans} and (\ref{eq:Mtrans}) is

\begin{corollary}\label{co:Mtrans}
 If $\eta$ and $\zeta$ are independent random variables, $\eta$ uniformly distributed on $(0,1)$ and $\zeta$  positive integer-valued on ${\mathbb N}$, then ${\mathfrak L}({\eta^{\frac{1}{\zeta}}}) = {\mathfrak L}(M[\zeta])$.
\end{corollary}

It is an {\em open problem}, if further transformations for nonnegative integer-valued random variables can be constructed to generate a random variable $\xi$ with $F_{\xi} = G$.

Another interesting problem is, whether an inverse of the $M$-transform exists in the following sense:

Let be given a random variable $\xi$ on $(0,1)$ with c.d.f. $F_{\xi}$. As in the construction of the $M$-transform, the i.i.d. random variables $\eta_1, \eta_2,...$ uniformly distributed on $(0,1)$ can be used as auxiliary quantities (or analogous to Corollary (\ref{co:Mtrans}) only one variable $\eta$). Now find  a function depending on $\xi$ and $\eta_1, \eta_2,...$ which defines a positive integer-valued $\zeta$ with generating function $G$ such that $G(x)=F_{\xi} (x)$ for all $x\in [0,1]$. From Corollary (\ref{co:Mtrans}) we obtain that both variables $\ln \xi$ and $\frac{1}{\zeta} \ln \eta$ have the same distribution. This is a mixture of exponential distributions with respect to $\zeta$.

\section{Nonnegative random variables and Laplace transform}

\subsection{Laplace transform and survival function}

For a nonnegative random variable $\zeta$ its Laplace transform $L_\zeta :[0,\infty )\to [0,1]$ is defined by
$$
L_\zeta (s)= {\mathbb E} {\rm e}^{-s \zeta} \quad \mbox{ for all } s\geq 0.
$$
As the Laplace transform depends only on the distribution of a random variable, we can also speak of the Laplace transform of a distribution.
Since $1-L_\zeta$ is continuous, nondecreasing with $1-L_\zeta (0)=0$ and $\lim_{s\to \infty} 1-L_\zeta (s)=1$, it can also be interpreted as the c.d.f. $F_\xi$ of some nonnegative random variable $\xi$, i.e. $F_\xi =1-L_\zeta$. Equivalently, $L_\zeta =1-F_\xi$ is the survival function of $\xi$.

{\em An open problem:} Let be given a nonnegative random variable $\zeta$ and a sequence $\eta_1, \eta_2, \ldots $ of i.i.d. random variables, uniformly distributed on $(0,1)$.
Find a random variable (if it exists) $\xi =\xi (\zeta  ,\eta_1, \eta_2, \ldots )$ which transforms $\zeta  ,\eta_1, \eta_2, \ldots $ such that $F_\xi =1-L_\zeta$. 
And as for the $M$-transform in  Section \ref{sec:M} we could ask for an inverse transform: For a given $\xi$ find a random variable $\zeta$ with Laplace transform $L_\zeta$ equal to the survival function of $\xi$.

\subsection{Laplace transform and probability generating function}

Recall that for a nonnegative integer-valued random variable $\zeta$  with generating function $G$, the Laplace transform is $L_\zeta (s)= G({\rm e}^{-s})$ for all $s\geq 0$.

Now we consider an arbitrary nonnegative random variable.

\begin{proposition}
Let $\zeta$ be a nonnegative random variable with Laplace transform $L_\zeta$ and define for all $t>0$ the function $G_t :[0,1] \to [0,1]$
by
$$
G_t(x)= L_\zeta (t(1-x)) \quad \mbox{ for all } x\in [0,1]. 
$$
\begin{enumerate}
\item Then for all $t>0$ the function $G_t$
is the generating function of a nonnegative integer-valued random variable.
\item If, for all $t>0$, $\kappa_t$ is a nonnegative integer-valued random variable with generating function $G_t$, then 
$$
\lim_{t\to \infty } L_{(\kappa_t/t)} (s)  = L_\zeta (s),
$$
which implies that for $t\to \infty$ the random variables $\kappa_t /t$  converge in distribution to $\zeta$.
\end{enumerate}

\end{proposition}

Proof:

As it can be seen in the proof of Theorem \ref{th:unbthick}, $G_t$ is the probability generating function of a nonnegative integer random variable, $\kappa_t$ say.
Now  define the nonnegative random variable $\beta_t =\kappa_t /t$ which has the Laplace transform $L_{\beta_t}$ with values
$$
L_{\beta_t} (s)= L_{\kappa_t} \left( \frac{s}{t} \right) = G_t \left( {\rm e}^{-\frac{s}{t}} \right) = L_\zeta \left( t\left( 1-{\rm e}^{-\frac{s}{t}}\right) \right) .
$$
This yields 
$$
\lim_{t\to \infty } L_{\beta_t} (s) = \lim_{a\to 0 } 
L_\zeta \left( \frac{1-{\rm e}^{-sa}}{a}  \right) = L_\zeta (s).
$$
\hfill $\Box$

{\it An open problem} is again the construction of the random variables $\kappa_t$ as a transform of a given $\zeta$.

\subsection{Roots of survival functions}

Let $\eta_1, \eta_2, \ldots $ be a sequence of i.i.d. nonnegative random variables with c.d.f. $F$. As it is well-known, for $n\in {\mathbb N} $ the survival function of the random variable
$\zeta :=\min \{ \eta_1,\ldots \eta_n\}$ is
$
1-F_\zeta (x)= (1-F(x))^n $ for all $x\geq 0 .
$ 
This immediately yields for the survival function of $\eta_1$ that
\begin{equation}\label{eq:survmin}
1-F = \sqrt[n] {1-F_\zeta} .
\end{equation}
 How can a random variable $\eta$ with c.d.f. $F$ according to (\ref{eq:survmin}) be generated from a sequence
$\zeta_1, \zeta_2, \ldots $ of i.i.d. copies of $\zeta$?

\begin{proposition}
Let $\zeta_1, \zeta_2, \ldots $ be a sequence of i.i.d. nonnegative random variables with c.d.f. $F_\zeta$ and $\alpha$ a random variable, geometrically distributed with parameter $1/n$, $n\in {\mathbb N}$, and independent from the sequence. Further, define the sequence  $\xi_1,\xi_2,\ldots $ of record times by 
\begin{eqnarray*}
\xi_1 & = & 1, \\
\xi_2 & = & \min \{ k> \xi_1: \zeta_k \geq \zeta_{\xi_1} \},\\
 &&\ldots \\
\xi_{m+1} & = & \min \{ k> \xi_m: \zeta_k \geq \zeta_{\xi_m} \},\\
&&\ldots
\end{eqnarray*}
Then the random variable $\zeta_{\xi_\alpha}$ has a c.d.f. $F$ satisfying (\ref{eq:survmin}).
\end{proposition}

Proof: 

As it is well-known (see e.g. \cite{resnick}), the process $\zeta_{\xi_1}, \zeta_{\xi_2},\ldots $ of records can be represented as a Poisson point process on $[0,\infty)$.

Given  $F_\zeta$, define the measure $\mu$ on $[0,\infty )$ (with the Borel $\sigma$-algebra) by
\begin{equation}\label{eq:failmeas}
\exp (-\mu ([0,x))) = 1- F_\zeta (x) \quad \mbox{ for all } x>0.
\end{equation}
This measure can be interpreted as a failure measure for $\zeta$. If $F_\zeta$ has the density $f_\zeta$, then for $x>0$ with $F_\zeta (x)<1$, the failure rate of $\zeta$ is $\frac{\partial \mu ([0,x))}{\partial x}=\frac{f_\zeta (x)}{1- F_\zeta (x)}$.  Note that $\mu$ is not necessarily a Radon measure.

Now let $\Psi$ be a Poisson point process on the positive half-axis with intensity measure $\mu$, and denote the ordered sequence of its points by $\beta_1 \leq \beta_2 \leq \ldots $ This implies ${\mathfrak L}(\beta_m) ={\mathfrak L}(\zeta_{\xi_m}) $ for $m=1,2,\ldots$
Now consider the Poisson point process $\Psi'$ generated from $\Psi$ by independent thinning with the probability $1- (1/n)$ for  deleting a point from $\Psi$. Then $\Psi'$ has the intensity measure $\mu' = (1/n)\mu $. Therefore, according to (\ref{eq:failmeas}) its first point $\beta_1'$ (in the ordered point set) has the c.d.f. satisfying (\ref{eq:survmin}). Furthermore, if $\alpha$ is geometrically distributed and independent from all the other random variables, we obtain that
${\mathfrak L}(\beta_1') ={\mathfrak L}(\zeta_{\xi_\alpha}) $, which completes the proof.
\hfill $\Box$

\subsection{Relations between exponential distributions}

It is well-known,  that the minimum of finitely many independent and exponentially distributed random variables is exponentially distributed as well. Furthermore, the sum of $n$ i.i.d.  exponentially distributed random variables with parameter $\lambda >0$ has an Erlang distribution, which is a special gamma distribution with parameters $n$ and $\lambda$.
In order to study the sum of not necessarily identically distributed random variables, we consider now particular convolutions of exponential distributions.

Some of the results have an interpretation concerning Poisson point processes. The intervals between the points of a homogeneous Poisson point process on the real axis are i.i.d. exponentially distributed. 

Denote the exponential distribution with parameter $\lambda >0$ by $\mathrm E[\lambda ]$ and by $\mathrm E^{ \ast k}[\lambda ]$ its $k$-fold convolution, $k\in {\mathbb N}$.

\begin{theorem}\label{theorem:convsum}
For all $0<\lambda <\infty$ and $0<p<1$, 
\begin{equation}\label{eq:convsum}
\mathrm E[p \lambda ] = p \sum_{k=0} ^{\infty} (1-p)^k \mathrm E^{ \ast (k+1)}[\lambda ]   .
\end{equation}
\end{theorem}

Proof: 

The proof is easy, using the Laplace transform $L(s)= \lambda /(\lambda +s)$, $s\geq 0$, for the exponential distribution with parameter $\lambda >0$, and the fact that the Laplace transform of a $k$-fold convolution of a distribution is just the $k$-th power of the Laplace transform of the respective distribution.  $\Box$

This result has also an interesting interpretation in terms of Poisson point processes on the positive real axis.  Let $\Phi$ be a homogeneous Poisson point process on $(0,\infty )$ with intensity $\lambda$. Then the coordinate of the first point of $\Phi$ has the exponential distribution $\mathrm E[ \lambda ]$, and the coordinate of the $(k+1)$-st point has the distribution $\mathrm E^{ \ast (k+1)}[\lambda ]$. Now consider the independent thinning of $\Phi$ where the points are deleted with probability $1-p$. This yields an homogeneous Poisson point process with intensity $p\, \lambda$. Thus the coordinate of the first point of the thinned point process  has the distribution $\mathrm E[p \lambda ]$. The probability that this first point of the thinned process (i.e. the first point which survived the independent thinning procedure) is the $(k+1)$-st point of $\Phi$ is $p(1-p)^k$. This is expressed by (\ref{eq:convsum}).

Decomposing the summands in (\ref{eq:convsum}) for $k\geq 1$ as
$$(1-p)^k \mathrm E^{ \ast (k+1)}[\lambda ]=(1-p)(1-p)^{k-1}\mathrm E[ \lambda ]\ast \mathrm E^{ \ast (k)}[\lambda ]$$ straightforwardly yields:
\begin{corollary}
For all $0<\lambda <\infty$ and $0<p<1$ 
$$
p \mathrm E[\lambda ] + (1-p) (\mathrm E[\lambda ]\ast \mathrm E[p \lambda ] )= \mathrm E[p \lambda ] .
$$
\end{corollary}

Substituting $\lambda$ by $\lambda_2$ and $p$  by $\lambda_1/\lambda_2$ for  $0<\lambda_1 < \lambda_2 <\infty$, this immediately supplies: 

\begin{corollary}\label{cor:faltung}
For all $0<\lambda_1 < \lambda_2 <\infty$
$$
 \mathrm E[\lambda_1 ] = \frac{\lambda_1}{\lambda_2} \mathrm E[ \lambda_2 ] + \frac{\lambda_2- \lambda_1}{\lambda_2} (\mathrm E[\lambda_1 ]\ast \mathrm E[ \lambda_2 ] )  
$$
or equivalently,
$$
\mathrm E[\lambda_1 ]\ast \mathrm E[ \lambda_2 ] =  \frac{\lambda_2}{\lambda_2- \lambda_1} \mathrm E[ \lambda_1 ] - \frac{\lambda_1}{\lambda_2- \lambda_1} \mathrm E[\lambda_2 ] .
$$

\end{corollary}

Now, we formulate the main result of this section for the convolution of two exponentially distributed random variables. Similarly as in Theorem \ref{theorem:convsum} it  is given as a mixture of Erlang distributions. Note that the two exponential distributions have different parameters.

\begin{theorem}\label{th:convolser}
For all $0<\lambda_1 < \lambda_2 <\infty$ and $p=\displaystyle{\left(\frac{\lambda_2 -\lambda_1}{\lambda_2 + \lambda_1}\right)^2}$ 
\begin{equation}\label{eq:convolser}
\mathrm E[\lambda_1 ] \ast \mathrm E[\lambda_2]
= \left(1-p \right) \, \sum_{k=0}^\infty 
p^k \; \mathrm E^{ \ast 2(k+1)}[{\textstyle \frac{1}{2}}(\lambda_2+ \lambda_1) ] .
\end{equation}

\end{theorem}

Proof:  

Let $L$ denote the Laplace transform of the distribution on the right-hand side of (\ref{eq:convolser}). Then, for $s\geq 0$ and $p=\displaystyle{\left(\frac{\lambda_2 -\lambda_1}{\lambda_2 + \lambda_1}\right)^2}$ ,
\begin{eqnarray*}
L(s) & =  & \left(  1-  \left( \frac{\lambda_2 - \lambda_1}{\lambda_2 + \lambda_1}\right)^2  \right) \sum_{k=0}^\infty 
\left( \frac{\lambda_2 - \lambda_1}{\lambda_2 + \lambda_1}\right)^{2k} \left( \frac{{\textstyle \frac{1}{2}} (\lambda_2+ \lambda_1)}{{\textstyle \frac{1}{2}}(\lambda_2+ \lambda_1)+s}\right)^{2k+2} \\ & &\\
&=& \left(  1-  \left( \frac{\lambda_2 - \lambda_1}{\lambda_2 + \lambda_1}\right)^2  \right) \left( \frac{\lambda_2 + \lambda_1}{\lambda_2+ \lambda_1 + 2s}\right)^{2} \frac{1}{1-\frac{(\lambda_2-\lambda_1)^2}{(\lambda_2 + \lambda_1 +2s)^2}} \\ & &\\
&=& \frac{\lambda_1 }{\lambda_1 +s} \cdot \frac{\lambda_2}{\lambda_2+s} ,
\end{eqnarray*}
and the term in the last line is just the product of the Laplace transforms of $\mathrm E[\lambda_1 ]$ and  $ \mathrm E[\lambda_2]$.
\hfill $\Box$

Alternatively, the result in Theorem \ref{th:convolser} also follows from an iterated application of the equation given in the next corollary.

\begin{corollary}
For all $0<\lambda_1 < \lambda_2 <\infty$ and $p=\displaystyle{\left(\frac{\lambda_2 -\lambda_1}{\lambda_2 + \lambda_1}\right)^2}$ 
\begin{equation}\label{eq:convexpo}
\mathrm E[ \lambda_1 ] \ast \mathrm E[ \lambda_2 ]
= \mathrm E^{ \ast 2}[{\textstyle \frac{1}{2}}(\lambda_1 + \lambda_2) ]\ast \Big( \left( 1- p \right) \delta_0 + p \, (\mathrm E[\lambda_1] \ast \mathrm E[\lambda_2 ]) \Big) .
\end{equation}

\end{corollary}
Again the proof is straightforward using the Laplace transforms.

\section{Concluding remarks}
In Joseph Mecke's fragments almost no references are given. Therefore we cannot reconstruct and cite the sources which he probably used. Consequently, we do not claim priority concerning all details.  
We are indebted to Hans Zessin for his valuable comments and hints.

\end{document}